\documentclass[10pt,a4paper]{amsart}
\usepackage{amssymb}
\usepackage{amsmath}
\usepackage{amsfonts}

\newcommand{\F}{\mathbb{F}}

\newtheorem{main-dummy}{Main-Dummy}
\newtheorem{dummy}{Dummy}

\newtheorem{main-theorem}[main-dummy]{Theorem}

\newtheorem*{theorem*}{Theorem}

\newtheorem*{cor*}{Corollary}

\theoremstyle{definition}

\theoremstyle{remark}
\newtheorem{rem}[dummy]{Remark}

\begin{document}
\bibliographystyle{amsalpha}

\author{Sandro Mattarei}

\email{mattarei@science.unitn.it}

\urladdr{http://www-math.science.unitn.it/\~{ }mattarei/}

\address{Dipartimento di Matematica\\
  Universit\`a degli Studi di Trento\\
  via Sommarive 14\\
  I-38050 Povo (Trento)\\
  Italy}

\title[On a bound of Garc\'{\i}a and Voloch]
{On a bound of Garc\'{\i}a and Voloch
for the number of points of a Fermat curve over a prime field}

\begin{abstract}
In 1988 Garc\'{\i}a and Voloch proved the upper bound $4n^{4/3}(p-1)^{2/3}$ for the number of solutions
over a prime finite field $\F_p$ of the Fermat equation $x^n+y^n=a$, where $a\in\F_p^\ast$ and
$n\ge 2$ is a divisor of $p-1$ such that $(n-\frac{1}{2})^{4}\ge p-1$.
This is better than Weil's bound $p+1+(n-1)(n-2)\sqrt{p}$ in the stated range.
By refining Garc\'{\i}a and Voloch's proof we show that the constant $4$ in their bound
can be replaced by $3\cdot 2^{-2/3}$.
\end{abstract}

\date{10 August 2005}



\thanks{The  author  is grateful  to  Ministero dell'Istruzione, dell'Universit\`a  e
  della  Ricerca, Italy,  for  financial  support of the
  project ``Graded Lie algebras  and pro-$p$-groups of finite width''.}

\maketitle

\thispagestyle{empty}

Let $\F_q$ be the finite field of $q$ elements  and let $p$ be its characteristic.
Consider the {\em Fermat curve}
$ax^n+by^n=z^n$, expressed in homogeneous coordinates, where $n>1$ is an integer prime to $p$,
and $a,b\in\F_q^\ast$. 
A classical estimate on the number $N_n(a,b,q)$ of its projective $\F_q$-rational points is
$
|N_n(a,b,q)-q-1|\le (n-1)(n-2)\sqrt{q}.
$
This is originally due to Hasse and Davenport~\cite{DavHas} but is a special case of
{\em Weil's bound} for curves over finite fields.
In the special case of Fermat curves Weil's bound is easy to prove by means of Gauss and Jacobi sums,
as well as its generalisation to {\em diagonal equations} in several variables,
see~\cite{IR}, \cite{LN} or~\cite{Small}.
An alternative proof uses character theory of finite groups,
see~\cite[Section~26]{Feit} for the basic idea and~\cite{Mat:Fermat-character}
for a refinement.

Weil's upper bound for $N_n(a,b,q)$ is not optimal when $n$ (and with it the genus of the curve)
is relatively large with respect to $q$.
Better upper bounds in this situation were found by Garc\'{\i}a and Voloch, using methods from algebraic geometry.
According to~\cite[Corollary~1]{GarVol}, rewritten here after elementary calculations,
if $s$ is an integer such that $1\le s\le n-3$ and $sn\le p$, then
\begin{equation}\label{eq:GV-bounds}
N_n(a,b,q)\le
\frac{1}{4}\left(
s^2-s-2+16\frac{1}{s+3}
\right)n^2
+2\frac{n(q-1-d)}{s+3}+d,
\end{equation}
where $d$ is the number of $\F_q$-rational points of the curve with $xyz=0$.
Garc\'{\i}a and Voloch pointed out that
their bounds~\eqref{eq:GV-bounds} hold in more general circumstances where
the assumption $sn\le p$ may not be satisfied,
and described those circumstances in detail for the cases $s=1,2$.
However, the special case stated above, and with $q=p$, was sufficient to them
for an application to Waring's problem in $\F_p$.
By estimating the minimum of their bounds, for $1\le s\le n-3$ and $sn\le p$,
they obtained the following intermediate result in~\cite[Section~3]{GarVol}:
the number of solutions $(x,y)\in\F_p\times\F_p$ of $x^n+y^n=a$, for $p$ a prime, $a\in\F_p^\ast$ and
$n\ge 2$ a divisor of $p-1$ such that $(n-\frac{1}{2})^{4}\ge p-1$, is at most
$4n^{4/3}(p-1)^{2/3}$.
A version of this bound (but for the equation $x^n-y^n=a$) with an unspecified constant in place of $4$
was later proved by Heath-Brown and Konyagin using Stepanov's method;
this is the case $T=1$ of~\cite[Lemma~5]{H-BK}, but see also~\cite[Chapter~3]{KonShp} for a generalization.
We comment further on this bound in Remark~\ref{rem:best}.
Mitkin has recently shown in~\cite{Mitkin} through elementary means that Garc\'{\i}a and Voloch's
bound holds (for the equation $x^n-y^n=a$) with the constant $4$ replaced by $2$, for $n>2^{3/4}(p-1)^{1/4}$.
However, it is also apparent from Garc\'{\i}a and Voloch's proof that the coefficient $4$ in their bound
can be lowered by refining their argument. 
In this note we bring the coefficient in that bound down to its optimal value
subject to being a consequence of the collection of Garc\'{\i}a and Voloch's bounds~\eqref{eq:GV-bounds}, as follows.

\begin{cor*}
Let $p$ be a prime and $a,b\in\F_p^\ast$.
Let $n\ge 4$ a divisor of $p-1$ such that $n^4\ge 4(p-1)$.
Then
$N_n(a,b,p)<3\cdot 2^{-2/3}n^{4/3}(p-1)^{2/3}$.
\end{cor*}

Since $3\cdot 2^{-2/3}$ is slightly less than $1.88989$, the Corollary is a little stronger than
the result in~\cite{Mitkin}.
We will deduce this result from the following more precise bound.

\begin{theorem*}
Let $p$ be a prime and $a,b\in\F_p^\ast$.
Let $n\ge 4$ a divisor of $p-1$ such that $n^4-2n^3-3n^2-8n\ge 4(p-1)$.
Then
\[
N_n(a,b,p)<
n^2\left(3(k/2)^{2/3}-
\frac{7}{2}(k/2)^{1/3}+
\frac{25}{12}\right),
\]
where $k=(p-1)/n$.
\end{theorem*}

\begin{proof}[Prof of the Theorem]
Because of the assumption $s\le n-3$, each bounding function in~\eqref{eq:GV-bounds}
does not decrease by replacing $d$ with its minimum value $0$.
We comment on the effect of this simplification in Remark~\ref{rem:d}.
In terms of $k$ and dropping the dependency on $d$ as described, the collection of bounds~\eqref{eq:GV-bounds} reads
\begin{equation}\label{eq:bound_min}
N_n(a,b,p)/n^2\le
\min\{U_s(k): 1\le s\le n-3,\, s\le k\},
\end{equation}
where
\[
U_s(k)=\frac{s^2-s-2}{4}
+2\frac{k+2}{s+3}.
\]
Thus, the upper bound for $N_n(a,b,p)/n^2$ given by inequality~\eqref{eq:bound_min} is a piece-wise linear
function of $k$.
Computation shows that
$U_{s+1}(k)=U_s(k)$ when $k=k_s$, where
$k_s+2=s(s+3)(s+4)/4$.
Because $k_k\ge k$ we have $U_s(k)\ge U_k(k)$ for $s\ge k$,
and hence the condition $s\le k$ is actually immaterial
in evaluating the minimum at the right-hand side of inequality~\eqref{eq:bound_min}.
It also follows that the right-hand side of inequality~\eqref{eq:bound_min}
is independent of $n$ for
$k\le k_{n-3}=\frac{1}{4}(n-3)n(n+1)-2$,
which is equivalent to our assumption
$n^4-2n^3-3n^2-8n\ge 4(p-1)$.
Therefore, under this assumption the bound~\eqref{eq:bound_min} can be written as
$N_n(a,b,p)/n^2\le V(k)$, where
$V(k)=\min\{U_s(k): s\ge 1\}$.

It remains to find a convenient function $W(k)$ which bounds the piece-wise linear function $V(k)$ from above.
Since
$V(k_s)=U_s(k_s)=(3s^2+7s-2)/4$, any concave function $W(k)$ such that
$W(k_s)\ge(3s^2+7s-2)/4$ for all integers $s\ge 1$ will do.

Consider the function
$W_c(k)=3(k/2)^{2/3}-
\frac{7}{2}(k/2)^{1/3}+c$,
where $c$ is a constant.
We have
\[
W_c((s+7/3)^3/4)=(3s^2+7s+4c)/4,
\]
and
$W_c'(k)=(k/2)^{-1/3}-
\frac{7}{12}(k/2)^{-4/3}\le
(k/2)^{-1/3}$.
In particular,
$W_c'(k_s)\le 2/s$
because $k_s\ge s^3/4$.
Since $W_c(k)$ is a concave function we have
\begin{align*}
W_c(k_s)&\ge W_c((s+7/3)^3/4)-
((s+7/3)^3/4-k_s)W_c'(k_s)
\\&\ge
\frac{3s^2+7s+4c}{4}
-\left(\frac{13}{12}s+\frac{559}{108}\right)\frac{2}{s}
=
V(k_s)
+c-\frac{5}{3}-\frac{559}{54 s}.
\end{align*}
Thus, if $c>5/3$ then $W_c(k_s)\ge V(k_s)$ for all integers $s\ge 1$ except a finite number.
A calculation now shows that the smallest value of $c$ such that
$W_c(k_s)\ge V(k_s)$ for all $s\ge 1$ is
$c=6-3(13/2)^{2/3}+(7/2)(13/2)^{1/3}$. (Equality then occurs for $s=2$.)
Since the value of this expression is (close to and) slightly less than $25/12$, the conclusion follows.
\end{proof}

\begin{rem}
The argument in the proof of the theorem can be extended to show that
$W_{71/48}(k)\le V(k)< W_{25/12}(k)$ for all $k\ge 1$.
The lower function equals the first three terms of the asymptotic expansion, for $k\to\infty$,
of the envelope of the family of linear functions
$U_s(k)$, where $s\ge 1$ is viewed as a real parameter instead of integral.
It follows that the bound for $N_n(a,b,p)$ given in the Theorem exceeds by less than $29n^2/48$
the minimum of the collection of bounds~\eqref{eq:GV-bounds}.
\end{rem}

\begin{rem}\label{rem:d}
We briefly explain the effect of having disregarded $d$ in the proof of the Theorem.
Let $G$ be the set of $n$th powers in $\F_p^\ast$, that is, the subgroup of $\F_p^\ast$ of order $k=(p-1)/n$.
If $(x,y)$ is a solution of $ax^n+by^n=1$ with $xy=0$ then any pair obtained from that by multiplying
$x$ and $y$ by elements of $G$ is also a solution.
Consequently, $N_n(a,b,p)-d$ is a multiple of $n^2$, and we can write~\eqref{eq:GV-bounds} in the form
\[
\big(N_n(a,b,p)-d\big)/n^2\le \left[U_s(k)-\tfrac{2}{s+3}(d/n)\right],
\]
where $U_s(k)$ as in the proof of the Theorem and with the square brackets denoting the integral part.
The ratio $d/n$ can only assume the values $0,1,2,3$, because it equals how many of
$a$, $b$ and $-a/b$ belong to $G$ (counting repetitions).
The proof of the Theorem (and, specifically, the formula for $k_s$)
shows that the strongest of the bounds~\eqref{eq:GV-bounds} for a given value of $k$
occurs, roughly, for $s$ close to $2(k/2)^{1/3}$.
Accordingly, one can improve the bound given in the Theorem
by making it dependent on $d$, but this would affect at most the term
$\tfrac{7}{2}(k/2)^{1/3}$, and not the leading term of the bound.
\end{rem}

\begin{proof}[Proof of the Corollary]
We only need to explain how the weakened conclusion allows us to relax our hypothesis
$n^4-2n^3-3n^2-8n\ge 4(p-1)$
to the weaker assumption $n^4\ge 4(p-1)$, which is equivalent to $n^3\ge 4k$.
When the stronger assumption is not satisfied, that is, when $k>k_{n-3}=\frac{1}{4}(n^3-2n^2-3n-8)$,
the bound~\eqref{eq:bound_min} reads
$N_n(a,b,p)/n^2\le U_{n-3}(k)$.
Thus, it suffices to show that
$U_{n-3}(k)<3(k/2)^{2/3}$
for $k_{n-3}<k\le n^3/4$.
Viewing $n$ as fixed, and hence $p$ as a function of $k$,
the left-hand side of the desired inequality is a linear function of $k$,
while the right-hand side is a concave function.
Since we know from the Theorem that the inequality is satisfied for $k=k_{n-3}$,
it remains only to check that this is the case also for $k=n^3/4$.
Indeed, we have
\[
U_{n-3}(n^3/4)^3=
\left(\frac{3n^2-7n+10}{4}+\frac{4}{n}\right)^3
<
\left(\frac{3n^2}{4}\right)^3
=
3\left(\frac{(n^3/4)}{2}\right)^2
\]
for all $n\ge 4$.
\end{proof}

\begin{rem}\label{rem:best}
The results from~\cite{H-BK} and~\cite{Mitkin} quoted in the introductory comments
both give upper bounds for the number of solutions $(x,y)\in\F_p^\ast\times\F_p^\ast$ of $x^n-y^n=a$,
for $p$ a prime, $a\in\F_p^\ast$ and
$n$ a divisor of $p-1$.
In particular, the special case $T=1$ of~\cite[Lemma~5]{H-BK} implies that
there is a constant $c$ such that the number of solutions is at most
$cn^{4/3}(p-1)^{2/3}$ if $n^4\ge p-1$.
An acceptable value for $c$ which follows from their proof is $4/(\sqrt{3}-1)$.
Our attempts to improve on this constant by refining their estimates
could not attain values lower than $2^{5/3}$, which is larger than $3$.

Mitkin's result in~\cite{Mitkin} is that the number of solutions is at most
$2n^{4/3}(p-1)^{2/3}$ if $n^4>8(p-1)$.
Although his method is very different from that of~\cite{GarVol},
Mitkin also establishes a family of bounds for the number of solutions divided by $n^2$,
which are linear in $k=(p-1)/n$ (like those of Garc\'{\i}a and Voloch summarized in Equation~\eqref{eq:bound_min}),
and then concludes by selecting the best of those for a given value of $k$.
However, Mitkin's family of bounds depends on three parameters rather than one,
and it seems not possible to individually match them with those of Garc\'{\i}a and Voloch.
The constant $2$ in Mitkin's final bound appears to be the best which can be attained by his method;
in fact, the stated purpose of Lemma~1 in~\cite{Mitkin} is to prove the bound with the constant $2$
rather than just $2+\varepsilon$ for some $\varepsilon>0$.
\end{rem}

\bibliography{References}

\end{document}